\documentclass[a4paper,11pt]{article}
\usepackage{amsmath}
\usepackage{amssymb}  
\usepackage[utf8x]{inputenc}
\newtheorem{prop}{Proposition}

\newtheorem{lem}[prop]{Lemma}

\newtheorem{thm}[prop]{Theorem}
\newtheorem{cor}[prop]{Corollary}
\newtheorem{defin}[prop]{Definition}

\newtheorem{quest}[prop]{Question}
\newcommand*\N{\mathbb{N}}

\newcommand*\Z{\mathbb{Z}}
\newcommand*\C{\mathbb{C}}

\begin{document}
\title{\large{\textbf{A CLASSIFICATION OF 5-DIMENSIONAL\\ MANIFOLDS, SOULS OF CODIMENSION TWO AND NON-DIFFEOMORPHIC PAIRS}}}
\author{\small{SADEEB OTTENBURGER}}
\maketitle
\begin{abstract}
Let $T(\gamma)$ be the total space of the canonical line bundle $\gamma$ over $\C P¹$ and $r$ an integer which is greater than one
and coprime to six. We prove that $L_r^3\times T(\gamma)$ admits an infinite sequence 
of metrics of nonnegative sectional curvature with pairwise non-homeomorphic souls, where $L_r^3$ is the standard 3-dimensional lens space with
fundamental group isomorphic to $\Z/r.$ We classify the total spaces of $S^1$-fibre bundles over $S^2\times S^2$ with fundamental group isomorphic to $\Z/r$ up to
diffeomorphism and use these results to give examples of manifolds $N$ which admit two complete metrics of nonnegative sectional curvature
with souls $S$ and $S'$ of codimension two such that $S$ and $S'$ are diffeomorphic whereas the pairs $(N,S)$ and $(N,S')$ are not diffeomorphic. 
This solves a problem posed by I. Belegradek, S. Kwasik and R. Schultz.
\end{abstract}
\section {Introduction}
In \cite{Ot-11} we classified the non-simply connected total spaces of principal $S^1$-fibre bundles over $S^2\times S^2$ with finite fundamental group
isomorphic to $\Z/r$, where $r$ is coprime to six, up to simple and tangential homotopy equivalence. This classification was a main step in
the discovery of the first examples of manifolds which admit infinitely many complete metrics of nonnegative sectional curvature and pairwise non-homeomorphic souls of codimension three. 
In this work we sharpen this result by proving that there exist manifolds which admit infinitely many complete metrics of nonnegative sectional curvature and pairwise
non-homeomorphic souls of codimension two. 
\\\\
Let $a,b$ be integers and $L^{a,b}$ be the total space of the principal $S^1$-fibre bundle over $S^2\times S^2$ given by the first Chern class $ax+by$, where $x$ 
and $y$ are the
standard generators of $H^2(\cdot;\Z)$ of the first and the second factor of the base respectively. And furthermore we denote the set $\{L^{a,b}\vert (0,0)\neq(a,b)\in \Z^2\}$ by 
$\cal L$. 
\begin{thm} Let $r,q$ be integers such that $r$ is greater than one and coprime to six.
The total space of a complex line bundle over $L^{r,qr}$ with primitive first Chern class admits an infinite
 sequence of complete metrics of nonnegative sectional curvature with pairwise non-homeomorphic souls.
\end{thm}
Let $T(\gamma)$ be the total space of the canonical line bundle $\gamma$ over $\C P^1$. If we consider $L^3_r$ as the total space 
of the principal $S^1$-fibre bundle over $S^2$ with first Chern class $rx$ then we obtain the following theorem as a special case of Theorem 1.
\begin{thm}                                 
If $r$ is as in Theorem 1 then $L_r^3\times T(\gamma)$ admits an infinite sequence of complete metrics of nonnegative sectional curvature with pairwise non-homeomorphic souls.
\end{thm}
Let $N$ be a smooth manifold and $\mathfrak{R}^c_{sec\geq 0}(N)$ be the set of smooth and complete metrics on $N$ of nonnegative sectional curvature with 
topology of smooth convergence on compact subsets. The diffeomorphism group of $N$ acts on $\mathfrak{R}^c_{sec\geq 0}(N)$ via pullback. The orbit
space under this action equipped with the quotient topology is called the \textit{moduli space} of complete metrics on $N$ of nonnegative sectional curvature
and it is denoted by $\mathfrak{M}^c_{sec\geq 0}(N)$.\\
The normal bundles of the souls that appear in Theorem 2 have nontrivial rational 
Euler class. By a result of V. Kapovitch, A. Petrunin and W. Tuschmann \cite[Lemma 5.1]{KPT-05} we know that the souls of metrics which lie in the same path component of $\mathfrak{M}^c_{sec\geq 0}(L_r^3\times T(\gamma))$ have to be diffeomorphic. This fact and Theorem 2 imply the following
\begin{prop}
If $r$ is as in Theorem 1 then $\mathfrak{M}^c_{sec\geq 0}(L_r^3\times T(\gamma))$ has infinitely many components.
\end{prop}
In this work we establish a classification theorem for a class of smooth closed 5-dimensional manifolds.
\begin{thm}Let $r$ be as in Theorem 1. Furthermore let $N$ and $N'$ be oriented smooth closed 5-dimensional manifolds which are homotopy equivalent 
to manifolds in $\cal L$ and which have fundamental groups isomorphic to $\Z/r$. Then $N$ and $N'$ are oriented diffeomorphic if and only if
there exists a orientation preserving homotopy equivalence $h:N'\rightarrow N$ such that
\begin{itemize}
\item [(i)] $\rho(h_*(g),N)=\rho(g,N')\textrm{ for all }g\in\pi_1(N')\setminus \{0\}$ \textrm{ and}
\item [(ii)] $\Delta(N)\sim\Delta (N'),$
\end{itemize}
where the $\rho$-invariant is defined in \cite[p.589]{AS-68} and $\Delta(\cdot)$ is the Reidemeister torsion as defined in \cite[p.405]{M-66}.
\end{thm}
We apply Theorem 4 to the classification of the non-simply connected ma\-ni\-folds in $\cal L$, where the order of their fundamental groups is coprime to six, 
up to diffeomorphism (Theorem 16). Then we use this explicit classification result and the surgery theory which lies behind it to prove Theorem 5.
\\\\
Let $x$ be an element of $ H^2(L^{r,qr};\Z)$ and $N_x^{r,qr}$ be the total space of the complex line bundle over $L^{r,qr}$ with first Chern class $x$.
\begin{thm} Let $r,q$ be integers, where $r$ is as in Theorem 1.\\\\
(i) If $q$ is not zero but divisible by $r$ and there is a unit $s$ in $\Z/r$ such that $s^2=-1$ then there exists a primitive element $y\in H^2(L^{r,qr};\Z)$ such that $N^{r,qr}_y$ has the following property:
$N^{r,qr}_y$ admits complete metrics $g$ and $g'$ of nonnegative sectional curvature with souls $S$, $S'$ respectively such that 
$S$ and $S'$ are diffeomorphic whereas the pairs $(N^{r,qr}_y,S)$ and $(N^{r,qr}_y,S')$ aren't diffeomorphic.\\\\
(ii) Let $x\in  H^2(L^{r,qr};\Z)$ be primitive and $g$, $g'$ be complete metrics of $N^{r,qr}_x$ of nonnegative sectional
curvature with souls $S$ and $S'$ respectively such that $S\in\cal L$. If $\Z/r$ does not contain a non-trivial unit $s$ with $s^3=1$ and $q$ is zero or not divisible by $r$ then 
$(N^{r,qr}_x,S)$ and $(N^{r,qr}_x,S')$ are diffeomorphic if and only if $S$ and $S'$ are diffeomorphic.
\end{thm}
If $r=5$ then $2$ is a unit in $\Z/5$ with $2^2=-1$ thus Theorem 5(i) yields solutions to Problem 4.9 in \cite{BKS1-09}. The statement of Theorem 5(i) also
holds if we require that there exists a non-trivial unit $s$ in $\Z/r$ with $s^2=1$. It is worth mentioning that Theorem 5 (ii) especially holds
for complete metrics on $L^3_r\times T(\gamma)$ of nonnegative sectional curvature.
\\\\The structure of this work and some remarks:
\\\\
\textit{Souls of codimension two.} We use new results of K. Grove and W. Ziller \cite {GZ-11} to prove the existence of complete metrics of nonnegative sectional curvature
on total spaces of complex line bundles over $L^{a,b}$ given by a pri\-mi\-tive integral cohomology class of degree two. Furthermore we apply techniques from surgery theory to prove Theorem 1.
\\\\
\textit{Topological Classification.} We prove Theorem 4 and classify the non-simply connected total spaces of principal $S^1$-fibre bundles over $S^2\times S^2$ with 
finite fundamental group of odd order not divisible by three up to diffeomorphism. A topological structure theorem is of course interesting for itself but it's nicer 
if it leads to new geometric insights. We apply our topological classification result and its proof in the context of nonnegative sectional curvature and prove Theorem 5. 
\\\\
By using similar methods as we use for proving Theorem 4 one can generalize these topological classification results
to a classification of smooth closed 5-manifolds with $\pi_1$ finite cyclic of odd order not divisible by three, $\pi_2$ torsion-free, where $\pi_1$ acts trivially on $\pi_2$.
Examples of such manifolds are total spaces of principal $S^1$-fibre bundles over simply connected smooth closed 4-manifolds with $\pi_2$ torsion-free 
and rank$(\pi_2)\geq 1$, e.g. $\C P^2\#\pm\C P^2$ or lens space bundles over $S^2$ with fibre a 3-dimensional lens space.\\
If we choose $\pi_2$ to be trivial then we are in the realm of 5-dimensional fake lens spaces which have been classified up to homeomorphism by C.T.C. Wall \cite[Part 3, 14 E ] {W-99}.
Thus the classification results here may be viewed as a special extension of Wall's results.
\\
This classification could also be of interest in contact topology \cite {GT-01} or might play a role in Sasakian geometry \cite {BG-08}.
\\\\\\
\textbf{Acknowledgements}\\\\
I thank my supervisor Matthias Kreck and Diarmuid Crowley for helpful discussions. And I am very grateful to Igor Belegradek for the idea of Theorem 1 and
useful comments.
\section{Souls of codimension two}
\begin{lem}
The total space of any complex line bundle over $L^{a,b}\in\cal L$ admits a complete metric of nonnegative sectional curvature with soul $L^{a,b}$.
\end{lem}
\textbf{Proof.} We know by definition that $L^{a,b}$ is the total space of a principal $S^1$-fibre bundle over $S^2\times S^2$ and its Gysin sequence implies that
$\Pi^*_{p,q}:H^2(S^2\times S^2;\Z)\rightarrow H^2(L^{a,b};\Z)$ is surjective. Thus any principal $S^1$-fibre bundle over $L^{a,b}$ is the pullback of a principal $S^1$-fibre
bundle over $S^2\times S^2$. Hence the total space of a principal $S^1$-fibre bundle over $L^{a,b}$ is (canonically) diffeomorphic to a principal $T^2$-bundle over $S^2\times S^2$.
From \cite[Theorem 4.5]{GZ-11} it follows that the total space of any principal $S^1$-fibre bundle over $L^{a,b}$ admits a complete metric of 
nonnegative sectional curvature which is invariant under the $T^2$-action. Let $\alpha$ be a complex line bundle over $L^{a,b}$, $T(\alpha)$ its total space and let $\alpha'$ be the corresponding principal $S^1$-bundle over $L^{a,b}$. Furthermore
let $\phi$ be the principal $T^2$-bundle over $S^2\times S^2$ with total space $T(\phi)$, where the inclusion of the first $S^1$-factor of $T^2$ into $T^2$ induces a $S^1$-fibre
bundle which is $\alpha'$. Thus $T(\alpha)$ and $T(\phi)$ are (canonically) diffeomorphic.
The total space $T(\alpha)$ can be obtained in the following way:
\[T(\alpha)=\C \stackrel{S^1\times \{1\}}{\times} T(\phi),\]
where we take the standard $S^1$-action on $\C$. Thus if we choose the standard metric on $\C$ and the $T^2$-invariant metric on $T(\phi)$ from \cite {GZ-11} then it 
follows from O'Neill's formula for the sectional curvature for Riemannian submersions that $T(\alpha)$ admits a complete metric of nonnegative sectional curvature, where
the soul is the zero-section of $\alpha$.\hfill$\square$\\\\
From \cite[Theorem 1]{Ot-11} we know that there are infinitely many pairwise non-homeomorphic total spaces of $S^1$-bundles over $S^2\times S^2$ which lie in the same 
tangential homotopy type of $L^{r,qr}$. By Corollary 19 we may choose the homotopy equivalences between these manifolds such that their normal invariants are trivial.
Following the idea of the proof of \cite[Theorem 14.1]{BKS1-09} we show the following: If we choose a complex line bundle over $L^{r,qr}$ with primitive first Chern class 
then all the total spaces of the pullback bundles with respect to homotopy equivalences with trivial normal invariant are diffeomorphic.\\\\
We give a sketch of the proof of \cite[Theorem 14.1]{BKS1-09} applied to our situation:\\
Let $\alpha$ be a complex line bundle over $L^{r,qr}$. By $D(\alpha)$ we denote the total space of the associated disc bundle, by $p$ the disc bundle projection
and by $S(\alpha)$ the boundary of $D(\alpha)$. If $f:N\rightarrow L^{r,qr}$ is a homotopy equivalence then $f$ induces a homotopy equivalence 
$f(\alpha):D(f^*\alpha)\rightarrow D(\alpha)$. The normal invariants of $f$ and $f(\alpha)$ are related by the following formula:
\begin{eqnarray}\mathfrak{q}(f(\alpha))=p^*\mathfrak{q}(f)\end{eqnarray}
which is for example proved in \cite[Lemma 5.9]{BKS-09}. Hence if $f$ has trivial normal invariant then the same holds for $f(\alpha)$.
\\
A \textit{simple homotopy structure} on $D(\alpha)$ is a simple homotopy equivalence $f:N\rightarrow D(\alpha)$, where $N$ is a smooth manifold.
Another simple homotopy structure $f':N'\rightarrow D(\alpha)$ and $f$ are said to be equivalent if there exists a 
diffeomorphism $d:N'\rightarrow N$ such that $f\circ d$ and $f'$ are homotopic. The set of equivalence classes of 
simple homotopy structures on $D(\alpha)$ is called the \textit{simple structure set} of $D(\alpha)$ and it is denoted by $S^s(D(\alpha))$ and the class which is represented 
by $id:D(\alpha)\rightarrow D(\alpha)$ is called the \textit{base point} of $S^s(D(\alpha))$. 
In Wall's surgery exact sequence associated to $S^s(D(\alpha))$ there is the map 
\[\Delta:L^s_{8}(\pi_1(D(\alpha)),\pi_1(S(\alpha)))\rightarrow S^s(D(\alpha))\]
which has the property that the elements of the image of $\Delta$ are represented by simple homotopy equivalences with trivial normal invariant and vice versa 
any homotopy equivalence with trivial normal invariant represents an element of the image of $\Delta$.
Thus if $L^s_{8}(\pi_1(D(\alpha)),\pi_1(S(\alpha)))$ is trivial then $f(\alpha)$ would represent the base point of
$S^s(D(\alpha))$ which means that $D(\alpha)$ and $D(f^*\alpha)$ would be diffeomorphic and hence the same would be true for their interiors.
If $i_*:\pi_1(S(\alpha))\rightarrow \pi_1(D(\alpha))$ is an isomorphism then it would
follow from Wall's $\pi$-$\pi$-Theorem that $L^s_{8}(\pi_1(D(\alpha)),\pi_1(S(\alpha)))$ is trivial.
\begin{lem}
Let $\alpha$ be a complex line bundle over $L^{r,qr}$ with primitive first Chern class then $i_*:\pi_1(S(\alpha))\rightarrow \pi_1(D(\alpha))$ is an isomorphism.
\end{lem}
\textbf{Proof.} The second stage of the Postnikov tower of $L^{r,qr}$ is $L^\infty_r\times \C P^\infty$ (\cite[Lemma 10]{Ot-11}, where $L^\infty_r$ is 
the infinite dimensional lens space with fundamental group isomorphic to $\Z/r$.
$H^2(L^\infty_r\times \C P^\infty;\Z)\cong H^2(L_r^\infty;\Z)\oplus H^2(\C P^\infty;\Z)$ which is isomorphic to $\Z/r\oplus \Z$. Let $z$ be the 
standard generator of $H^2(\C P^\infty;\Z)$ and $v_1$ a generator of $H^1(L_r^\infty;\Z/r)$. Let $\left (H^1(L^{r,qr};\Z/r)\right )^*$ be the set of units of
$H^1(L^{r,qr};\Z/r)$ and $\mathrm{Pri}(H^2(L^{r,qr};\Z))$ the set of elements in $H^2(L^{r,qr};\Z)$ which generate a $\Z$-summand.
By \cite[Lemma 15]{Ot-11} there is a 1-1 correspondence between 
\[\left (H^1(L^{r,qr};\Z/r)\right )^*\oplus \mathrm{Pri}(H^2(L^{r,qr};\Z))\]
and the set of homotopy classed of
 maps $g$ from $L^{r,qr}$ to $L^\infty_r\times \C P^\infty$ which induce isomorphisms on $\pi_1$ and $\pi_2$. The correspondence is given by
\begin{eqnarray*}
g^*(v_1)\textrm{ and }g^*(z).
\end{eqnarray*}
Thus a complex line bundle over $L^{r,qr}$ with primitive first Chern class $\bar{z}$ is the pullback of a complex line bundle over $L^\infty_r\times \C P^\infty$ under a map $g$ 
with $g^*(z)=\bar{z}$. The homotopy exact sequence associated to the principal $S^1$-fibre bundle over $L^\infty_r\times \C P^\infty$ with first Chern class $z$ starts as
follows:
\begin{eqnarray*}
 \pi_2(\C P^\infty)\stackrel{\partial}{\rightarrow}\pi_1(S^1)\rightarrow \pi_1(S_z)\rightarrow \pi_1(L^\infty_r)\rightarrow 0,
\end{eqnarray*}
where $S_z$ is the associated sphere bundle. The boundary map $\partial$ is an isomorphism since it is induced by the Euler class of the bundle. Hence $\pi_1(S_z)\rightarrow \pi_1(L^\infty_r)$ is an isomorphism. 
But by using the 2-smoothing $g$ the analogues statement transfers to the $S^1$-fibre bundle over $L^{r,qr}$ given by the first Chern class $\bar{z}$.
\hfill$\square$
\section {Topological classification}
Let $L^{a,b}\in\cal L$. Before we give a proof of Theorem 4 we gather some differential topological properties of $L^{a,b}$ (\cite [Section 2]{Ot-11}):
\begin{itemize}
 \item [i)] $\pi_1(L^{a,b})\cong \Z /\textrm{gcd}(a,b)$ and $\pi_2(L^{a,b})\cong\Z$.
\item [ii)] The tangent bundle of $L^{a,b}$ is stably trivial, which implies that all its stable characteristic classes are trivial.
\item [iii)] The Reidemeister torsion as defined in \cite[p.405] {M-66} is trivial.
\end{itemize}
We work in the smooth category of closed 5-manifolds and within this ca\-te\-gory non-diffeomorphic always implies non-homeomorphic which relies on 
the fact that there are no exotic smooth structures on $S⁵$.
\\\\
How surgery gets involved in the classification of manifolds of the type we consider was partially explained in section 2 of \cite {Ot-11}. There
we found a connection between a bordism classification and a homotopy classification:
\\\\
Let $\pi_1(L^{a,b})\cong\Z/r$. The second Postnikov stage of $L^{a,b}$ is $L_r^{\infty}\times \C P^{\infty}=:B_r$ (\cite[Lemma 10] {Ot-11}).
The bordism group of our interest is $\Omega^{Spin}_5(B_r)$ which
is defined to be the set 
\[\left\{(M,f)\vert \begin{array}{cc}
	M\textrm{ a closed smooth} 
                            \textrm{ 5-dimensional spin manifold},\\
			      f:M\rightarrow B_r\end{array}\right \} \]
modulo an equivalence relation
which is given as follows: $(M,g)\sim(N,g')$ if there exists a 6-dimensional smooth spin manifold with boundary equal to the disjoint union of $M$ and $N$ and a map
 $F:W\rightarrow B_r$
which restricted to the boundary components is the map $f$, $f'$ respectively.\\\\
The Postnikov decomposition of $L^{a,b}$ yields maps $f:L^{a,b}\rightarrow B_r$ which induce isomorphisms on the first 
and second homotopy groups. We call such maps \textit{normal 2-smoothings}. Let $f:L^{a,b}\rightarrow B_r$ and $f':L^{a',b'}\rightarrow B_r$ be normal 2-smoothings
 then $(L^{a,b},f)$, $(L^{a',b'},f')$ define elements in $\Omega^{Spin}_5(B_r)$. Assume
that $(L^{a,b},f)$ and $(L^{a',b'},f')$ represent the same element in $\Omega^{Spin}_5(B_r)$ then there exists a bordism $(W,F)$ between them as
described above and we call $(W,F)$ a \textit{normal (co)bordism}. Furthermore we may assume that $F$ is a 3-equivalence \cite [Prop. 4]{K-99}. 
\\\\
Let $\Lambda$ be the group ring $\Z[\Z/r]$ and let $K\pi_3(W)$ be $\mathrm{Ker}(F_*:\pi_3(W)\rightarrow \pi_3(B_r)),$ $K\pi_3(L^{a,b})$ be $\mathrm{Ker}(f_*:\pi_3(L^{a,b})\rightarrow \pi_3(B_r))$ 
and $\mathrm{Im}K\pi_3(L^{a,b})$ be the image of $K\pi_3(N_i)$ under the homomorphism which is induced by the inclusion $L^{a,b}\hookrightarrow W.$
In \cite [Section 2]{Ot-11} we explained how we associate to $(W,F)$ an element $\Theta (W,F)$ in a group called $L^{s,\tau}_6(\Z/r,S)$:
\begin{eqnarray*}
\left[\bar{\lambda}:\frac{K\pi_3(W)}{\mathrm{Im}K\pi_3(N_0)}\times \frac{K\pi_3(W)}{\mathrm{Im}K\pi_3(N_0)}\rightarrow \Lambda, \tilde{\mu}\right]=:\Theta(W,F).
\end{eqnarray*}
If there exists another normal bordism $(W',F')$ between $(L^{a,b},f)$ and $(L^{a',b'},$ $f')$ 
then $\Theta(W,F)=\Theta(W',F')$ and $(W,F)$ and $(W',F')$ differ by a sequence of surgeries re\-la\-tive the boundary $(L^{a,b},f)\cup(L^{a',b'},f')$. 
\begin{lem}
Up to normal $s$-cobordism
an element $\Theta$ in $L^{s,\tau}_6(\Z/r,S\oplus \Z)$ uniquely determines a normal 2-smoothing $(N,g)$, i.e. if $(V,G)$ is a normal bordism $(V,G)$
between $(L^{a,b},f)$ and $(N,g)$ then $\Theta(V,G)=\Theta$. 
\end{lem}
The proof of all this can be found in \cite [Theorem 5.8]{W-99} where it is done for normal $5$-smoothings 
but it extends literally to our situation.\\\\
The group $L^{s,\tau}_6(\Z/r,S\oplus \Z)$ is related to the ordinary Wall group $L^{s}_6(\Z/r,S\oplus \Z)$ by the following exact sequence:
\[0\rightarrow L^s_6(\Z/r,S\oplus \Z) \stackrel{i}{\rightarrow} L^{s,\tau}_6(\Z/r,S\oplus \Z)\stackrel{\tau}{\rightarrow}\textrm{Wh }(\Z/r),\]
where $i$ is just the canonical inclusion and the map $\tau$ sends stable equivalence classes of weakly based non-singular $(-1)$-quadratic forms 
over $\Lambda$ to the Whitehead torsion of the matrix representation of the adjoint of this quadratic form with respect to the chosen weak equivalence class of basis.
\\
More details to all the technical definitions may be found in \cite{K-99} or \cite[Section 2]{Ot-11}.
\begin{thm}$(W,F)$ is bordant relative boundary to an $s$-cobordism if and only if $\Theta (W,F)=0.$
\end{thm}
\textbf{Proof.} 
The proof of this theorem is analogues to the proof of Theorem 3 in \cite {K-99}.
\hfill$\square$
\\\\
The so called \textit{multisignature} extends the notion of the ordinary signature of a quadratic form over the integers in a certain sense. The most general definition of the
multisignature which one can find in \cite {W-66} or \cite [p.174]{W-99} is applicable to unimodular forms over group rings of finite groups with either the standard or a
 non-standard involution. Since we deal with forms over group rings with trivial orientation ("all manifolds are orientable") we give a definition of the multisignature which
 is equivalent to the general definition restricted to the orientable case \cite [p.175] {W-99}:\\\\
Let ${\cal M}$ be a free $\Lambda $-module and $\lambda:{\cal M}\times {\cal M} \rightarrow \Lambda$ a non-degenerate skew-hermitian form. 
Furthermore we denote by $p_c$ the map from $\Lambda$ to $\Z$ which sends $a_0+\sum_{\alpha\in \Z/r \setminus\{0\}}a_\alpha$ to $a_0.$ The composition $p_c \circ \lambda$ is 
a skew-hermitian $\Z$-valued  non-degenerate form. We extend $p_c \circ \lambda$ to ${\cal M}^{\C}:={\cal M}\otimes_{\Z}\C$ in the obvious way which yields a 
skew-hermitian unimodular $\C$-valued form $\lambda_{\C}$, i.e. 
$\lambda_{\C} (x,y)=-\overline{\lambda_{\C} (y,x)}$ for all $x,y\in {\cal M}^{\C}.$\\
It is clear that ${\cal M}^{\C}$ inherits a $\Z/r$-action from ${\cal M}$ and we easily realize that
$\lambda_{\C}(x\alpha,y\alpha)=\lambda_{\C}(x,y),$ for all $x,y\in {\cal M}^{\C}\textrm{ and }\alpha\in \Z/r.$
Now we choose a positive definite $\Z/r$-invariant hermitian form $\left\langle   \cdot ,  \cdot \right\rangle_{\Z/r}$ on ${\cal M}^{\C}$:
Let $(\tilde{z}_1,...,\tilde{z}_l)$ be a complex basis of ${\cal M}^\C$ and let $\left\langle \cdot,\cdot\right\rangle $ be the standard hermitian form on
 ${\cal M}^\C$, i.e. $\left\langle \tilde{z}_i,\tilde{z}_j\right\rangle =\delta_{ij}.$ We define the following $\Z/r$-invariant hermitian product:
$\left\langle \cdot,\cdot\right\rangle _{\Z/r}:=\sum_{\alpha\in \Z/r}\left\langle \alpha(\cdot),\alpha(\cdot)\right\rangle.$ 
There is the following linear map $A$ of ${\cal M}^{\C}$ to itself defined by $\lambda_{\C} (x,y)=\left\langle x,Ay\right\rangle\,$ for all $x,y\in {\cal M}^{\C}.$ 
It follows that $A$ is a skew-hermitian $\Z/r$-equivariant automorphism of ${\cal M}^{\C}$ and that the eigenvalues of $A$ are purely imaginary and nonzero.
 Thus ${\cal M}^{\C}={\cal M}_+^{\C}\oplus {\cal M}_-^{\C},$ where ${\cal M}_{\pm}^{\C}$ is the sum of the eigenspaces corresponding to positive multiples of
 $\pm i.$ Since the eigenspaces ${\cal M}_{+}^{\C}$ and ${\cal M}_{-}^{\C}$ are $\Z/r$-invariant they define complex $\Z/r$-re\-pre\-sen\-ta\-tions. We denote
the characters of these $\Z/r$-representations by $r_+$ and $r_-$ respectively.
\begin{defin}
1) The multisignature of $\lambda$ which we denote by $MS(\lambda)$ is the element of the complex representation ring $RU(\Z/r)$ given by 
\[r_+  - r_-.\]
2) Let $g$ be an element of $\Z/r$ then $MS(g,\lambda):=r_+(g)  -r_-(g).$
\end{defin}
The multisignature is well defined. This follows from the following three facts:
\begin{itemize}
\item [i)] The characters depend continuously on the inner product.
\item [ii)] The space of all $\Z/r$-invariant hermitian products on ${\cal M}^\C$ is connected.
\item [iii)] The characters of a compact group are discrete.
\end{itemize}
Let $r$ be odd. If $\Z/r$ is the trivial group one
knows that the signature induces an isomorphism between $L_{4m}^s(\{1\})$ and $8\Z.$ It is also true that the so called Arf-invariant gives an isomorphism
between $L_{4m+2}^s(\{1\})$ and $\Z/2.$ The multisignature is an important tool for distinguishing elements in $L_{2k}^s(\Z/r):$
\begin{thm}(\cite [Thm. 13A.4.] {W-99})There is a decomposition 
\[L_{2k}^s(\Z/r)=L_{2k}^s(\{1\})\oplus \widetilde {L}_{2k}^s(\Z/r),\]
where the \textit{multisignature} maps $\widetilde {L}_{2k}^s(\pi)$ injectively to the characters (real or imaginary as appropriate) trivial on 1 and divisible by 4.
\end{thm}
In the following we study the surgery obstruction group $L^{s,\tau}_{6}(\Z/r,S\oplus \Z)$.
\\\\ 
First we have a closer look on $L^{s}_{6}(\Z/r)$:\\
As already indicated in Theorem 11 the functorial character of the Wall $L$-groups yields a decomposition of
$L^s_{6}(\Z/r)$, i.e. $L^s_{6}(\Z/r)=L^s_{6}(\{1\})\oplus \widetilde {L}^s_{6}(\Z/r).$ 
Let $[(\Lambda^d,\lambda,\mu)]=:E$ be an element in $L^s_{6}(\Z/r).$ The first coordinate of $E$ with respect to the above decomposition is
detected by the Arf-invariant in the following sense:\\
To an element $[(\Lambda^d,\lambda,\mu)]$ in $L^s_{6}(\Z/r)$ one can assign an element in $L^s_{6}(\{1\})$:\\ As we identify $\Lambda$ with $\Z^b$ (for some $b\in\N$) in a 
canonical way we regard $\Lambda^d$ as a $\Z$-module. Let $\epsilon:\Lambda\rightarrow \Z$ be the augmentation map which is a ring homomorphism and let
 $\tilde{\epsilon}$ be the projection of $\frac{\Lambda}{\langle2a\rangle}$ to $\frac{\epsilon(\Lambda)}{\epsilon(\langle2a\rangle)}\cong \Z/2.$\\
We compose $\lambda$ with $\epsilon$ and compose the quadratic refinement $\mu$ with $\tilde{\epsilon}$ then we obtain an integral non-degenerate skew-hermitian
quadratic form which represents an element $e$ in $L^s_{6}(\{1\})$. We define Arf$(E)$ to be the \textit{classical Arf-invariant} of $e$. The second coordinate of $E$ with respect to
the decomposition above is according to Theorem 11 detected by the multisignature.\\
The difference between $L^s_{6}(\Z/r)$ and $L^{s}_{6}(\Z/r,S\oplus \Z)$ comes from the different choices of form parameters. In order to understand 
$L^{s}_{6}(\Z/r,S(W)\oplus \Z)$ we do the same considerations as for $L^s_{6}(\Z/r).$ We observe that because of the $\Z$ in $S(W)\oplus \Z$ elements in
$L^{s}_{6}(\Z/r,S(W)\oplus \Z)$ are uniquely determined by the multisignature since the Arf-invariant of elements in $L^{s}_{6}(\{1\},S(W)\oplus \Z)$ is trivial.\\
If the multisignature of $(\bar{\lambda}:\Lambda^d\times \Lambda^d\rightarrow \Lambda,\bar{\mu})$ is zero, then this implies that there is a basis ${\cal B}$ of $\Lambda^d$ such that
$(\bar{\lambda},\bar{\mu})$ represents the zero-element in $L^{s}_{6}(\Z/r,S(W)\oplus \Z)$ and hence in $L^{s,\tau}_{6}(\Z/r,S(W)\oplus \Z)$ if one chooses $\Lambda^d$ with the
basis $\cal B$. 
Let $(\bar{\lambda},\bar{\mu})$ be the equivariant intersection form and the self-intersection associated to $(W,F)$ which represents the 
surgery obstruction $\Theta(W,F)$ (\cite [Section 2]{Ot-11}).
The notations $(\bar{\lambda},\bar{\mu})$ and $(\bar{\lambda},\bar{\mu})_{{\cal B}}$ shall indicate that $\Lambda^d$ is equipped with a preferred basis and the basis ${\cal B}$ 
respectively.
\\
By abuse of notation a canonically based $(-1)$-hyperbolic form over $(\Lambda,S\oplus \Z)$ is again denoted by $H^n_{-1}(\Lambda)$ (for some $n\in\N$).
If $(\bar{\lambda},\bar{\mu})_{{\cal B}}$ represents the zero element in $L^{s}_{6}(\Z/r,S\oplus \Z)$, then there is a $l\in\N$ such that
$H^l_{-1}(\Lambda)\oplus (\bar{\lambda},\bar{\mu})_{{\cal B}}$ is isomorphic to $H^{l'}_{-1}(\Lambda)$ for some $l'\in \N$, where the isomorphism is simple.
Let $B$ be the matrix of base change with respect to the bases in $H^l_{-1}(\Lambda)\oplus (\bar{\lambda}, \bar{\mu})$ and 
$H^l_{-1}(\Lambda)\oplus (\bar{\lambda}, \bar{\mu})_{{\cal B}}.$ The element in Wh$(\Z/r)$ that is represented by $B$ is denoted by $\tau(B)$ and we define $\mathrm{MS}(\Theta(W,F))$ 
to be $\mathrm{MS}(\bar{\lambda}, \bar{\mu})$ and $\mathrm{MS}(\alpha,\Theta(W,F))$ to be $\mathrm{MS}(\alpha,\bar{\lambda}, \bar{\mu}).$
These considerations imply
\begin{prop}
$\Theta (W,F)=0$ if and only if $\mathrm{MS}(\Theta (W,F)) \textrm{ and }\tau (B)$ are trivial.
\end{prop}
The proof of Theorem 3 in \cite {K-99} implies the following
\begin{thm}$W$ is relative boundary bordant to an $h$-cobordism $(W_h;L^{a,b},$ $L^{a',b'})$ if and only if the multisignature of the equivariant intersection pairing $\bar{\lambda}$
on $\frac{\pi_3(W)}{\mathrm{Im}K\pi_3(L^{a,b})}$ is trivial. In this case the vanishing of the algebraic torsion $\tau (B)$ is equivalent to the vanishing of the Whitehead torsion
of the inclusion $L^{a,b}\hookrightarrow W_h.$
\end{thm}
Before we begin with the proof of Theorem 4 we explain what the \textit{equi\-va\-riant $\Z/r$-signature} of a $6$-dimensional bordism $W$ on which $\Z/r$
operates smoothly and in an orientation preserving way is. For more details we refer to \cite {AS-68}.
\\\\
We know that there is the following skew-hermitian form called the intersection form of $W$:
\[\lambda\ :\ H^3(W,\partial W;\Z)\times H^3(W,\partial W;\Z)\rightarrow \Z.\]
This form is $\Z/r$-invariant, i.e. $\lambda(gu,gv)=\lambda(u,v)$. This follows from the assumption that $\Z/r$ acts on $W$ in an orientation preserving way.\\
The radical $\mathrm{rad}(\lambda)$ of $\lambda$ equals $\mathrm{Ker}(i^* : H^3(W,\partial W;\Z)\rightarrow H^3(W,\Z))$. We easily see that the form 
\[
\bar{\lambda}\ : \ \frac{H^3(W,\partial W;\Z)}{\mathrm{rad}(\lambda)}\times\frac{H^3(W,\partial W;\Z)}{\mathrm{rad}(\lambda)}\rightarrow \Z,
\]given by
\[\bar{\lambda}([u],[v]):=\lambda(u,v)\]
is a non-degenerate $\Z/r$-invariant skew-hermitian form.\\ 
In the following we denote $\frac{H^3(W,\partial W;\Z)}{\mathrm{rad}(\lambda)}$ by $\hat{H}(W)$.
Tensoring $\hat{H}(W)$ with $\C$ over $\Z$ yields a complex vector space which we denote by $\hat{H}(W)^\C$. We extend $\bar{\lambda}$ to $\hat{H}(W)^\C$ in the following way:
\[
\bar{\lambda}_\C\ :\ \hat{H}(W)^\C\times \hat{H}(W)^\C\rightarrow\C,\ (x\otimes z_1,y\otimes z_2)\mapsto {\bar{\lambda}}(x,y)\cdot(z_1\cdot z_2).
\]
This complex valued quadratic form is a skew-hermitian $\Z/r$-invariant unimodular form.
\\
As in the definition of the multisignature we obtain unitary $\Z/r$-representations $r_i$, $r_{-i}$ which lead us to
\begin{defin}
The $\Z/r$-\textbf{signature} $\mathrm{sign}(\Z/r,W)$ of $W$  is defined as
\begin{eqnarray*}
r_i-r_{-i}\in RU(\Z/r)
\end{eqnarray*}
Let $g\in \Z/r$ then 
\[\mathrm{sign}(g,W):=r_i(g)-r_{-i}(g).\]
\end{defin}
\textbf{Proof of Theorem 4.}
\begin{prop}(\cite[Prop.12,Cor.13]{Ot-11})
Let $r$ be as in Theorem 1 and $N,M$ be closed smooth oriented spin 5-manifolds with vanishing first Pontrjagin classes and $f$ and $g$ normal
 2-smoothings from $N$ and $M$ to $B_r$ respectively. Then the following statements are equivalent:
\begin{itemize}
\item [i)] $(N,f)$ and $(M,g)$ represent the same element in $\Omega^{Spin}_5(B_r)$.
\item [ii)] $f_*[N]=g_*[M]\in H_5(B_r;\Z)$.
\item [ii)] $N$ and $M$ are homotopy equivalent.
\end{itemize}
\end{prop}
Let us assume that there exists a normal bordism $(W,F)$ between normal 2-smoothings $(N,f)$ and $(N',f')$ and we may choose $(W,F)$ such that $F$ is a 3-equivalence
(see \cite[Prop.15]{K-99}).
We identify $\pi_1(W)$ with $\Z/r$ and denote $\Z[\Z/r]$ by $\Lambda.$ In the following we relate the multisignature of $\Theta(W,F)$ to the $\Z/r$-equivariant 
signature of $\widetilde{W}$. We equip $W$ and $\widetilde {W}$ with basepoints $b$ and $\tilde {b}$ respectively such that $\tilde{b}$ is a lift of $b$ under the universal
 covering map. We want to study the unimodular skew-hermitian form
\[\bar{\lambda} :\frac{K\pi_3 (W)}{\mathrm{Im} K\pi_3 (N)}\times \frac{K\pi_3 (W)}{\mathrm{Im} K\pi_3 (N)}\rightarrow \Lambda\]
which comes from the intersection pairing defined on $K\pi_3(W).$ We observe that $K\pi_3(W):=\mathrm{Ker}(F_{\star}:\pi_3(W)\rightarrow \pi_3(B_r))\stackrel{\pi_3 (B_r)=0}{=}{\pi_3(W)}.$ 
\\
There is the following composition of maps
\begin{eqnarray*}
&\pi_3(\widetilde{W})\stackrel{{\cal H}}{\rightarrow}H_3(\widetilde{W};\Z)\stackrel{{(\cap [\widetilde{W},\partial \widetilde{W}])^{-1}}}
{\rightarrow} H^3(\widetilde{W},\partial \widetilde{W};\Z)&\rightarrow
\\&\rightarrow\underbrace{\frac{H^3(\widetilde{W},\partial \widetilde{W};\Z)}
{\mathrm{Ker}(i^*:H^3(\widetilde{W},\partial \widetilde{W};\Z)\rightarrow H^3(\widetilde{W};\Z))}}_{=:\hat{H}^3(\widetilde{W})}&\end{eqnarray*}
which we call $\Phi,$ where $(\cap [\widetilde{W},\partial \widetilde{W}])^{-1}$ is the inverse of the Poincar\'e-Lefschetz isomorphism.
The map $\Phi$ is $\Z/r$-equivariant and since $\widetilde{W}$ is 1-connected we know that the Hurewicz map ${\cal H}$ is surjective.
We have seen that the $\Z/r$-signature for $\widetilde {W}$ is defined for the non-singular pairing
$\bar{{\cal I}} :\hat {H}^3(\widetilde{W})\times\hat {H}^3(\widetilde{W})\rightarrow \Z$
which comes from the cup-product-pairing 
on $H^3(\widetilde{W},\partial \widetilde{W};\Z)$. We denote the intersection pairing on $K\pi_3 (W)=\pi_3 (W)$ by $\lambda.$ Let us recall how $\lambda$ was defined.
For a detailed exposition of the following we refer to \cite [Ch.5] {W-99}. \\
Elements of $\pi_3 (W)$ are regular homotopy classes of immersions of $S^3$ into $W.$
For $\gamma_0$, $\gamma_1\in  \pi_3 (W)$ we find representatives $(i_0,w_0)$ and $(i_1,w_1)$ such that the images of these maps intersect transversally only in
 finitely many points. We denote the set of intersection points by $D$. From $\gamma_0$ and $\gamma_1$ we obtain a well defined element in $\Lambda$, namely 
$\lambda(\gamma_0,\gamma_1):=\sum_{d\in D}\epsilon (d)\alpha(d),$ where $\epsilon (d)\in\{\pm 1\}$ and $\alpha(d)\in \Z/r$ (\cite [p. 45] {W-99}). There is a unique
lift $\widetilde{i_j}$ of $i_j$ to $\widetilde {W}$ determined by $\tilde {b}.$ Let $\lambda_{\Z}([\widetilde{i_0}],[\widetilde{i_1}])$ be the $\Z$-valued algebraic
 intersection number of the homology classes $[\widetilde{i_0}]$ and $[\widetilde{i_1}]$ then we may identify $\lambda(\gamma_0,\gamma_1)$ with
 $\sum _{\alpha\in \Z/r}\lambda_{\Z}([\widetilde{i_0}],[l_{\alpha^{-1}}\circ\widetilde{i_1}])\alpha,$
where $l_{\alpha^{-1}}$ denotes the left multiplication by $\alpha^{-1}.$ This means that
\[p_c \circ \lambda (i_0,i_1)=\lambda_{\Z}([\widetilde{i_0}],[\widetilde{i_1}]).\]
But $\lambda_{\Z}([\widetilde{i_0}],[\widetilde{i_1}])=[\widetilde{i_0}]^*\cup[\widetilde{i_1}]^*\in H^6(\widetilde{W},\partial \widetilde{W};\Z),$
where $[\widetilde{i_j}]^*$ denotes the Poincar\'e-Lefschetz dual of $[\widetilde{i_j}].$ All in all we obtain
\[p_c\circ \lambda=\bar{{\cal I}}(\Phi(\cdot),\Phi(\cdot)).\]
Thus the radicals of $p_c\circ\lambda$, $\lambda$ and $\bar{{\cal J}}(\Phi(\cdot),\Phi(\cdot))$ coincide and are equal to $\mathrm{Ker}(\Phi)$ which is $\mathrm{Im}K\pi_3(N).$ From the
 definitions of the multisignature and the equivariant signature we conclude:\\\\
\textit{Computing the multisignature of $\bar{\lambda}$ is the same as computing the $\Z/r$-sig\-na\-ture of $\widetilde{W}.$}\\\\
Thus
\begin{eqnarray}MS(\alpha, \Theta (W,F))=\mathrm{sign}(\alpha,\widetilde{W}),\ \ \forall \alpha\in \Z/r.\end{eqnarray}
Let $\alpha\in\Z/r\setminus\{0\}$ then from (2) and Novikov's additivity formula for the equivariant signature it follows that 
\[MS(\alpha, \Theta (W,F))=\rho (\alpha,N)-\rho (\alpha,N').\]
Since $\widetilde{W}$ is 6-dimensional it follows that $\mathrm{sign}(\widetilde{W})$ is trivial and thus
\[MS(\alpha, \Theta (W,F))=0,\ \ \forall \alpha\in \Z/r\]
if and only if 
\begin{eqnarray}\rho (\beta,N')=\rho ((f_*^{-1}\circ f'_*)(\beta),N), 
\forall \beta\in \pi_1(N')\setminus \{0\}.\end{eqnarray}
If the $\rho$-invariant condition holds then we know from Theorem 13 that $W$ is normally bordant to an $h$-cobordism.\\
Since $\pi_1(N)$ and $\pi_1(N')$ are cyclic. By assumption $N$ and $N'$ are homotopy equivalent to manifolds in $\cal L$ thus $\pi_1(N)$ and $\pi_1(N')$
operate trivially on the rational cohomology ring of the universal covering spaces $\widetilde{N}$ and $\widetilde{N'}$ respectively. These facts imply 
that the Reidemeister torsions $\Delta(N)$ and $\Delta(N')$ according to \cite[p.405]{M-66} are defined.
Since by assumption $\Delta(N)\sim\Delta(N')$ \cite[Thm. 12.8]{M-66} imply that the Whitehead torsion of the inclusion of one boundary component into this $h$-cobordism is trivial. By applying the $s$-cobordism theorem
$N$ and $N'$ are diffeomorphic.
\\
The theorem follows by the following considerations: If there is a homotopy equivalence $h:N'\rightarrow N$ and $f:N\rightarrow B_r$ 
is a normal 2-smoothing then $f':=f\circ h: N'\rightarrow B_r$ is a normal 2-smoothing of $N'$. Condition (i) in Theorem 4 is equivalent to condition (3). 
And by the considerations above $N$ and $N'$ are diffeomorphic if conditions (i) and (ii) of Theorem 4 hold.
\hfill$\square$\\\\
Our next aim is to classify the non-simply connected manifolds in $\cal L$, where the order of the fundamental group is coprime to 6. The orientation convention for
the manifolds in $\cal L$ is the following: Since these manifolds are total spaces of principal $S^1$-fibre bundles over $S^2\times S^2$ we orient them by orienting the base
and the fibre in the standard way.
\begin{thm}
Let $r$ be as in Theorem 1 and $L^{a,b}$, $L^{a',b'}\in{\cal L}$ with $\pi_1(L^{a,b})\cong\pi_1(L^{a',b'})\cong \Z/r$. Furthermore let $(m,n)$, $(m',n')$
 be pairs of integers such that $m\frac{b}{r}+n\frac{a}{r}=1=m'\frac{b'}{r}+n'\frac{a'}{r}$. Then 
$L^{a,b}$ and $L^{a',b'}$ are diffeomorphic if and only if there exist $\epsilon,\epsilon ',\delta\in \{\pm 1\}$ and $k,k'\in \Z/r$ such that
\[ab=\delta a'b',\ \ \ \ \ \ \ \ \ \ \ \ \ \ \ \ \ \ \ \]
 \[(\epsilon m + k\frac{a}{r})(\epsilon n - k \frac{b}{r})\equiv  \delta (\epsilon' m' + k'\frac{a'}{r})(\epsilon' n' - k' \frac{b'}{r})\mathrm{mod\textrm{ }} r,\ \ \ \ \ \ \ \]
 \[\frac{b}{r}(\epsilon m +k\frac{a}{r})-\frac{a}{r}(\epsilon n -k\frac{b}{r})\equiv
 (\frac{b'}{r}(\epsilon' m' +k'\frac{a'}{r})-\frac{a'}{r}(\epsilon n' -k'\frac{b'}{r}))\mathrm{mod\textrm{ }}r.\]
 \end{thm}
\textbf{Proof of Theorem 16.}
\begin{thm}(\cite[Theorem 8]{Ot-11}) Let $r$ be as in Theorem 1 and $L^{a,b}$, $L^{a',b'}\in{\cal L}$ with $\pi_1(L^{a,b})\cong\pi_1(L^{a',b'})\cong \Z/r$.
 Furthermore let $(m,n)$, $(m',n')$ be pairs of integers such that $m\frac{b}{r}+n\frac{a}{r}=1=m'\frac{b'}{r}+n'\frac{a'}{r}$.
 Then $L^{a,b}$ and $L^{a',b'}$ are oriented homotopy equivalent if and only if there exist $s,s'\in(\Z/r)^*$, $\epsilon,\epsilon '\in \{\pm 1\}$ and $k,k'\in \Z/r$ such that
\[s^3\frac{ab}{r^2}\equiv  s'^3\frac{a'b'}{r^2}\mathrm{mod\textrm{ }}r,\ \ \ \ \]
\[s (\epsilon m + k\frac{a}{r})(\epsilon n - k \frac{b}{r})\equiv  s' (\epsilon' m' + k'\frac{a'}{r})(\epsilon' n' - k' \frac{b'}{r})\mathrm{mod\textrm{ }}r,\]
\[s^2(\frac{b}{r}(\epsilon m +k\frac{a}{r})-\frac{a}{r}(\epsilon n -k\frac{b}{r}))\equiv  s'^2(\frac{b'}{r}(\epsilon' m' +k'\frac{a'}{r})-\frac{a'}{r}
(\epsilon n' -k'\frac{b'}{r}))\mathrm{mod\textrm{ }}r.\]
\end{thm}
If we change the orientation of one of the manifolds in Theorem 16 and Theorem 17 then we have to insert a minus sign on the corresponding side of the congruences.
Since it is important for later discussions we explain where the parameters $\epsilon, k$ and $s$ come from:\\\\
Let us fix a choice of tuples $(m,n)$ as in Theorem 16. By \cite [Lemma 15]{Ot-11} we can explicitly give a 1-1 correspondence between the set of homotopy classes of normal 2-smoothings of $L^{a,b}$
and the set of triples
\[\{(\epsilon, k,s)\vert \epsilon \in\{\pm 1\}, k\in\Z/r,s\in(\Z/r)^*\}\]
by sending a normal 2-smoothing $g:L^{a,b}\rightarrow B_r$ to 
\begin{eqnarray}
&&g^*(v_1)=s\alpha,\\
&& g^*(z)=\epsilon(m\Pi_{a,b}^*(x)-n\Pi_{a,b} ^* (y))+k(\frac{a}{r}\Pi_{a,b}^*(x)+\frac{b}{r}\Pi_{a,b}^*(y)),
\end{eqnarray}
where $\Pi_{a,b}$ is the projection map of the associated $S^1$-fibre bundle of $L^{a,b}$ and $v_1$ and $z$ are generators of $H¹(L^\infty_r;\Z/r)$ and
$H^2(\C P^{\infty};\Z)$ respectively.
\\\\
We know that the universal covering space of $L^{a,b}$ is $L^{\frac{a}{r}\frac{b}{r}}$,
and that the deck transformation on $D^{\frac{a}{r}\frac{b}{r}}$ by $\pi_1(L^{a,b})$ is given by fibrewise rotation by angles corresponding to the $r$'th roots of unity. This 
perspective yields a canonical identification of $\pi_1(L^{a,b})$ with $\Z/r$. The disc bundle $D^{\frac{a}{r}\frac{b}{r}}$ associated to the $S^1$-fibre bundle structure
with the $\Z/r$-action canonically extended serves as a convinient choice of a bordism. Furthermore this $\Z/r$-bordism has trivial equivariant signature since
on the one hand the $\Z/r$-action is homotopically trivial, as it sits in an $S^1$-action and on the other hand the dimension 
of the bordism is not divisible by 4 which means that the ordinary signature is trivial. The fixed point set is just $S^2\times S^2$ and
the normal bundle of the fixed point set is isomorphic to the 2-dimensional real vector bundle given by the Euler class $\frac{a}{r}x+\frac{b}{r}y$.
 Let $g$ be a non-trivial element of $\Z/r$ and $\theta_g$ the rotation angle between $0$ and $\pi$ of the action by $g$ then the $\rho$-invariant associated 
to the action of $g$ is defined to be the evaluation of certain characteristic polynomials depending on the Chern-, Pontrjagin classes of 
the normal bundle $\cal N$ of $S^2\times S^2$ and the Pontrjagin classes of $S^2\times S^2$, on the fundamental class of $S^2\times S^2$:
\begin{eqnarray*}
(i \tan \frac{\theta _\beta}{2})^{-1}\sum _{j=0}{\cal L}_j(S^2\times S^2)
\sum _r {\cal M}^{\theta _\beta}_r ({\cal N}_{\theta _\beta})\left[ S^2\times S^2\right ]_{\pm},
\end{eqnarray*}
where $\left[ S^2\times S^2\right ]_{\pm}$ is $\pm 1$ times the standard fundamental class of $S^2\times S^2.$ The sign in $\left [ S^2\times S^2\right ] _{\pm}$ 
depends on how $\beta$ operates on $L^{\frac{a}{r},\frac{b}{r}}$ and we conclude that
\begin{eqnarray}\rho (-g,L^{\frac{a}{r},\frac{b}{r}})=-\rho (g,L^{\frac{a}{r},\frac{b}{r}}).\end{eqnarray}
The formula for the $\rho$-invariant is the following:
\begin{eqnarray}\rho:\pi_1(L^{a,b})\setminus \{0\}\rightarrow \C,\ \ g\mapsto-i\frac{\cos (\frac{\theta_g}{2})}{2r^2\sin^3 (\frac{\theta_g}{2})}ab\end{eqnarray}
(\cite [p.88]{Ot-09})
which is clearly injective if $ab\neq 0$. If the $\rho$-invariant of $L^{a,b}$ and $L^{a',b'}$ shall coincide then the non-triviality of (7) 
implies that $|ab|=|a'b'|$. And (6) together with the injectivity of (7) implies that $s'=\delta s$ if $ab=\delta a'b'$, where $\delta\in\{\pm 1\}.$
Hence the conditions of Theorem 4 are fulfilled for $L^{a,b}$ and $L^{a',b'}$ if and only if there exists $m,m',n,n',k,k',\epsilon,\epsilon'$ and $\delta\in\{\pm 1\}$ 
as stated in Theorem 16 such that
\[ \frac{ab}{r^2}=\delta\frac{a'b'}{r^2},\ \ \ \ \ \ \ \ \ \ \ \ \ \ \ \ \ \ \ \]
 \[(\epsilon m + k\frac{a}{r})(\epsilon n - k \frac{b}{r})\equiv  \delta (\epsilon' m' + k'\frac{a'}{r})(\epsilon' n' - k' \frac{b'}{r})\mathrm{mod\textrm{ }} r,\ \ \ \ \ \ \ \]
 \[\frac{b}{r}(\epsilon m +k\frac{a}{r})-\frac{a}{r}(\epsilon n -k\frac{b}{r})\equiv \underbrace{\delta^2}_{=1}
 (\frac{b'}{r}(\epsilon' m' +k'\frac{a'}{r})-\frac{a'}{r}(\epsilon n' -k'\frac{b'}{r}))\mathrm{mod\textrm{ }}r.\]
\hfill$\square$\\\\
\textbf{Proof of Theorem 5 (i).} The idea of the proof of Theorem 5 (i) is the following: We prove for $L^{r,qr}$ that there is a self-homotopy equivalence $h$
with trivial normal invariant and a primitive element $y$ of $H²(L^{r,qr};\Z)$ such that there is no self-diffeomorphism $\phi$ such that 
\begin{eqnarray}h^*(y)=\pm\phi^*(y).\end{eqnarray} 
Since $h$ has trivial normal invariant the total spaces of the complex line bundles with first Chern class $y$ and $h^*(y)$ respectively are diffeomorphic (proof of Theorem 1).
But if there exists a
self-diffeomorphism of the total space which sends the zero-section of the one bundle to the zero-section of the other then $\phi$ would preserve the normal bundle structure, 
i.e. formula (8) would hold which can not be the case.\\\\
Let $(W,F)$ be a normal bordism between $(L^{r,qr},f_0)$ and $(L^{r,qr},f_1)$, where $f_0,f_1:L^{r,qr}\rightarrow B_r$ are normal 2-smoothings. We may assume that $F$ is 
a 3-equivalence. As we have already seen the surgery obstruction $\Theta (W,F)$ is determined by its multisignature and may be canonically identified with an
element in $L^s(\Z/r,S)$. As already mentioned to any
 element $\Theta'$ of $L^{s,\tau}(\Z/r,S)$ there exists a normal bordism $(W',F')$ between $(L^{r,qr},f_0)$ and $(N,f')$, where $f'$ is a normal 2-smoothing
and $(N,f')$ is unique up to normal $s$-cobordism. Let us take the element $\Theta_s$ of $L^s(\Z/r,S)$ which has the same multisignature (up to isomorphism) as $\Theta (W,F)$.
Then Wall's surgery implies: There exists
\begin{itemize}
\item [i)] a bordism  $(W_s;L^{r,qr},N_s)$,
\item [ii)] a map $F_s:W_s\rightarrow L^{r,qr}\times I$, 
where $F_s\vert_{L^{r,qr}}$ is the identity and $F_s\vert_{N_s}$ is a simple homotopy equivalence and
\item [iii)] a stable normal framing $\Phi_s$ of $F_s ^* \nu_{L^{r,qr}}\oplus \tau_{W_s}$,
where $\nu_{L^{r,qr}}$ is the stable normal bundle of $L^{r,qr}$.
\end{itemize}
$F_s\vert_{N_s}$ is unique up to normal $s$-cobordism. And thus $\Theta_s$ maps to an element in $S^s(L^{r,qr})$, the simple structure set of $L^{r,qr}$.
Let $F_0:L^{r,qr}\times I\rightarrow B_r$ be a homotopy of $f_0$ then $F_0\circ F_s:W_s\rightarrow B_r$ is normal bordism between
 $(L^{r,qr},f_0)$ and $(N_s,F_0\circ F_s\vert_{N_s})$. These considerations together with Lemma 8 yield the following
\begin{lem}
 $(N_1,f_1)$ and $(N_s,F_0\circ F_s\vert_{N_s})$ differ by a normal $s$-cobordism thus the map $f_1$ can be obtained from $f_0$ by precomposing with a simple self-homotopy equivalence 
with trivial normal invariant.
\end{lem}
Let us have a look on Theorem 17 in the case where $L^{a,b}=L^{a',b'}=L^{r,qr}$ and we choose $m=0$ and $n=1$.
Recall that the further parameters which are used in Theorem 16 determine the images of generators of $H¹(L^\infty_r;\Z/r)$ and $H²(\C P^\infty;\Z)$ under the map which is induced
by a normal 2-smoothing ((4),(5)) and if we choose the opposite orientation we have to insert a minus to the
corresponding sides of the congruences.
We claim that if $q$ is not zero but divisible by $r$ then the congruences in Theorem 17 hold for: $s$, $s'=1$, $\epsilon=1$, $\epsilon'=-1$, $k=-1$ and $k'=s$: The
first equation is trivially fulfilled since $q\equiv 0\mathrm{\textrm{ }mod\textrm{ }}r$. The other two congruences are:
\begin{eqnarray}
 sk\epsilon\equiv s'k'\epsilon'\textrm{mod }r,\\
s^2\epsilon \equiv s'^2\epsilon'\textrm{mod }r.
\end{eqnarray}
We realize that these congruences hold for the above choice of parameters. Hence by Lemma 18 there is a self-homotopy equivalence $h$ with trivial normal invariant such that
\begin{eqnarray*}
 h^*(m\Pi_{r,rq}^*(x)-n\Pi_{r,rq}^*(y)-\Pi_{r,rq}^*(x)+q\Pi_{r,rq}^*(y))
\end{eqnarray*}
equals
\begin{eqnarray*}
 -(m\Pi_{r,rq}^*(x)-n\Pi_{r,rq}^*(y))+s(\Pi_{r,rq}^*(x)+q\Pi_{r,rq}^*(y))
\end{eqnarray*}
which obviously is not $\pm 1$ times its preimage.
But if there was a self-diffeomorphism $\phi$ fulfilling the same cohomological property as $h$ then we could construct a normal $s$-cobordism:\\\\
Glue $W_1:=(L^{r,qr}\times I_1;L^{r,qr}\times \{0\},L^{r,qr}\times \{1\})$
and $W_2:=(L^{r,qr}\times I_2;L^{r,qr}\times \{0\},L^{r,qr}\times \{1\})$ together by identifying $L^{r,qr}\times \{1\}$ with $L^{r,qr}\times \{0\}$ via $\phi$. We call the result
$W_1\cup_{\phi}W_2$ and by $D$ we denote the canonical diffeomorphism from $W_1\cup_{\phi}W_2$ to $L^{r,qr}\times I$. Let $f:L^{r,qr}\rightarrow B_r$ be the normal
 2-smoothing which corresponds to the choice of the parameters $s,\epsilon$ and $k$ ((3),(4)) and $F:L^{r,qr}\times I\rightarrow B_r$ be a homotopy of $f$. 
Then $(W_1\cup_{\phi}W_2, F\circ D)$ is a normal $s$-cobordism
between $(L^{r,qr},f)$ and $(L^{r,qr},F\circ D\vert_{L^{r,qr}\times \{1\}})$, where $F\circ D\vert_{L^{r,qr}\times \{1\}}$ is the normal 2-smoothing from $L^{r,qr}$ to $B_r$
which corresponds to the choice of the parameters $s',\epsilon'$ and $k'$. But since $r^2q\neq0$ the $\rho$-invariant is injective and since $(W_1\cup_{\phi}W_2, F\circ \delta)$ is a normal $s$-cobordism the 
 the equations in Theorem 17 should hold for $s=\pm s'$ and $\epsilon',k'$ as above (proof of Thm. 16). 
Equation (9) implies that $s=k'=\pm k=\pm 1$ but $s^2= - 1$ which is a contradiction.\\\\
If there is a non-trivial unit $s$ in $\Z/r$ with $s^2=1$ then we choose the parameters as follows:$s$, $s'=1$, $\epsilon=1=\epsilon'$, $k=1$ and $k'=s$. 
They fulfill the congruences (9) and (10) and the statement of Theorem 5(i) follows from the considerations above.
 \hfill $\square$
\\\\
\textbf{Proof of Theorem 5 (ii).} If $(N^{r,qr}_x,S)$ and $(N^{r,qr}_x,S')$ are diffeomorphic then trivially $S,S'$ are diffeomorphic.\\\\
From \cite[Cor. 5.12]{BKS-09} it follows that there is a homotopy equivalence $h:S\rightarrow S'$ which pulls back the normal bundles.
\\
If $q$ is zero then we realize that since the $\rho$-invariant is zero any isomorphism on the second cohomology group induced
by a self-homotopy equivalence can be induced by a self-diffeomorphism. This follows from the proof of Theorem 16 and Lemma 18, where the simple homotopy equivalence is a diffeomorphism.\\
If $q$ is not divisible by three and there is no non-trivial unit $s$ in $\Z/r$ such that $s^3=1$ then the first equation in Theorem 17 implies 
that $s'=\pm s$ depending on the chosen orientations. 
But the above considerations imply that any set of parameters with $s'=\pm s$ fulfilling the equations in Theorem 17 also fulfills the equations in Theorem 16. Hence again
any isomorphism on the second cohomology group induced by a self-homotopy equivalence can be induced by a self-diffeomorphism.
\hfill $\square$
\begin{cor} Let $r$ be as in Theorem 1.
 All the manifolds in $\cal L$ which are homotopy equivalent to $L^{a,b}$ with $\pi_1(L^{a,b})\cong\Z/r$ lie in the orbit of the basepoint of the simple
 structure set of $L^{a,b}$ under the action of $L^s_6(\Z/r)$. Thus the homotopy equivalences between manifolds in $\cal {L}$ may be chosen in 
such a way that they have trivial normal invariant.
\end{cor}
\textbf{Proof.} If $L^{a,b}$ and $L^{a',b'}$ are homotopy equivalent then there exists a normal bordism $(W,F)$ between normal 
2-smoothings of $L^{a,b}$ and $L^{a',b'}$. The considerations before Lemma 8 and Lemma 18 imply that there exists a simple homotopy equivalence between 
$L^{a,b}$ and $L^{a',b'}$.\hfill$\square$\\\\
The discussions and the results in this work motivate some interesting questions which relate the simple homotopy type to curvature.\\\\
Let $M$ and $N$ be smooth closed non-simply connected manifolds with finite cyclic fundamental group isomorphic to $\pi$.
\begin{quest}Is it true that
if $M$ and $N$ differ by an element of $L_6^s(\pi)$ and $M$ admits a metric of nonnegative sectional curvature then $N$ admits a metric 
of nonnegative sectional curvature?
\end{quest}
More generally one could ask the following
\begin{quest}
Is it true that if $M$ and $N$ are simply homotopy equivalent and $M$ admits a metric of nonnegative sectional curvature then $N$ 
admits a metric of nonnegative sectional curvature?
\end{quest}


\begin{thebibliography}{[]}

\bibitem[AS-68]{AS-68} M.F. Atiyah, I.M. Singer,  \textit{The index of elliptic operators: III}, Ann. of Math. 87, 546-604, 1968.
\bibitem [BG-08]{BG-08} C. Boyer and K. Galicki, \textit{Sasakian geometry}, Oxford University Press, 2008.
\bibitem[BKS-09]{BKS-09} I. Belegradek, S. Kwasik, R. Schultz, \textit{Moduli spaces of nonnegative sectional curvature and non-unique souls}, arXiv:0912.4869 (math.DG), 2009.
\bibitem[BKS1-09]{BKS1-09} I. Belegradek, S. Kwasik, R. Schultz, \textit{Codimension two souls and cancellation phenomena}, arXiv:0912.4874v1 (math.DG), 2009.
\bibitem[GT-01]{GT-01} H. Geiges and C. B. Thomas, \textit{Contact structures, equivariant spin bordism, and periodic fundamental groups}, Math. Annalen 320, 685-708, 2001.
\bibitem[GZ-11]{GZ-11} K. Grove, W. Ziller, \textit{Lifting group actions and nonnegative curvature}, Trans. Amer. Math. Soc., 2011.
\bibitem[K-99]{K-99} M. Kreck, \textit{Surgery and duality}, Ann. of Math. 149, 707-754, 1999.
\bibitem[KPT-05]{KPT-05} V. Kapovitch, A. Petrunin, W. Tuschmann, \textit{Nonnegative pinching, moduli
spaces and bundles with infinitely many souls}, J. Differential Geom. 71, no. 3, 365-383, 2005.
\bibitem[M-66]{M-66} J. Milnor, \textit{Whitehead torsion}, Bull. Amer. Math. Soc. 72, 385-426, 1966.
\bibitem[Ot-09]{Ot-09} S. Ottenburger, \textit{A diffeomorphism classification of 5- and 7-dimensional non-simply-connected homogeneous spaces}, http://hss.ulb.uni-bonn.de/2009/1820/1820.pdf.
\bibitem[Ot-11]{Ot-11} S. Ottenburger, \textit{Simply and tangentially homotopy equivalent but non-homeomorphic homogeneous manifolds}, arXive preprint, 2011.
\bibitem[Sh-79]{Sh-79} V. A. Sharafutdinov, \textit{Convex sets in a manifold of nonnegative curvature}, Mat. Zametki
26, no. 1, 129-136, 159, 1979.
\bibitem[W-99]{W-99} C.T.C. Wall, \textit{Surgery on compact manifolds} 2nd Edition (edited by A.A. Ranicki), Mathematical Surveys and Monographs Vol 69, 1999.
\bibitem[W-66]{W-66} C.T.C. Wall, \textit{Surgery of non-simply connected manifolds}, Ann. of Math. 84, 217-276, 1966.

\end{thebibliography}
\end{document}